\documentstyle[12pt]{article}

\newtheorem{thm}{Theorem}

\newtheorem{pr}[thm]{Proposition}

\newenvironment{pf}{\noindent {\em Proof:}}{$\Box$\\}

\newcommand{\N}{\mbox{\hskip.1em N \hskip -1.25em \relax I \hskip .1em}}

\newcommand{\summ}{\sum^{\infty}_{n=m}}
\newcommand{\sumum}{\sum^{m-1}_{n=1}}
\newcommand{\sumnj}{\sum^k_{j=1}t_{n_j}e_j}
\newcommand{\sumi}{\sum^i_{j=1}t_{n_j}e_j}
\newcommand{\sumMnj}{\sum^k_{j=1}M(t_{n_j})}
\newcommand{\sumj}{\sum^k_{j=1}}
\newcommand{\lan}{\langle}
\newcommand{\ran}{\rangle}
\newcommand{\ep}{\epsilon}
\newcommand{\cz}{c_0}
\newcommand{\pn}{p.n.\ }
\newcommand{\ipn}{i.p.n.\ }
\newcommand{\iso}{isomorphically }
\newcommand{\poly}{polyhedral }
\newcommand{\llb}{\biggl\{}
\newcommand{\lrb}{\biggr\}}
\newcommand{\fateq}{\hspace{1em} = \hspace{1em}}
\newcommand{\fatleq}{\hspace{1em} \leq \hspace{1em}}
\newcommand{\fatarr}{\hspace{1em} \to \hspace{1em}}

\newcommand{\spn}{\mbox{span}}

\begin{document}

\begin{center}
  {\Large\bf Some isomorphically polyhedral\\
             Orlicz sequence spaces}\vspace{3mm}\\
  {\large\sc Denny H. Leung}
\end{center}

\vspace{1mm}

\begin{abstract}
A Banach space is polyhedral if the unit ball of each of its
finite dimensional subspaces is a polyhedron.  It is known that a
polyhedral Banach space has a separable dual and is $c_0$-saturated,
i.e., each closed infinite dimensional subspace contains an isomorph
of $c_0$.  In this paper, we show that the Orlicz sequence space $h_M$
is isomorphic to a polyhedral Banach space if $\lim_{t\to 0}M(Kt)/M(t)
= \infty$ for some $K < \infty$.
We also construct an Orlicz sequence space $h_M$ which is
$c_0$-saturated, but which is not isomorphic to any polyhedral Banach
space.  This shows that being $c_0$-saturated and having a separable
dual are not sufficient for a Banach space to be isomorphic
to a polyhedral Banach space.
\end{abstract}

\begin{figure}[b]
  \rule{3in}{.005in}
\end{figure}

A Banach space is said to be {\em polyhedral} if the
unit ball of each
of its finite dimensional subspaces is a polyhedron.  It is {\em
isomorphically polyhedral} if it is isomorphic to a polyhedral Banach
space.  Fundamental results concerning polyhedral Banach spaces were
obtained by Fonf \cite{F1, F2}.

\begin{thm}\label{c0sat}
{\rm (Fonf)} A separable isomorphically polyhedral Banach space
is $c_0$-saturated and has a separable dual.
\end{thm}

Recall that a Banach space is $\cz$-{\em saturated} if every closed
infinite dimensional subspace contains an isomorph of $\cz$.
Fonf also proved a characterization of isomorphically polyhedral
spaces in terms of certain norming subsets in the dual.  In order to
state the relevant results, we introduce some terminology due to
Rosenthal \cite{R1, R2}. The (closed) unit ball of a Banach space $E$
is denoted by $U_E$.\\

\noindent{\bf Definition\ } Let $E$ be a Banach space.\\
(1) A subset $W \subseteq E'$ is {\em precisely norming} (p.n.) if
$W \subseteq U_{E'}$, and for all $x \in E$, there is a $w \in W$ such
that $\|x\| = |w(x)|$.\\
(2) A subset $W \subseteq E'$ is {\em isomorphically precisely
norming} (i.p.n.) if $W$ is bounded and \\
(a) there exists $K < \infty$ such
that $\|x\| \leq K\sup_{w \in W}|w(x)|$ for all $x
\in E$, \\
(b) the supremum $\sup_{w \in W}|w(x)|$ is attained
at some $w_0 \in W$ for all $x \in E$.\\

It is easy to see that $W \subseteq E'$ is \ipn if and only if there
is an equivalent norm $|||\cdot|||$ on $E$ so that $W$ is \pn in
$(E,|||\cdot|||)'$.

\begin{thm}\label{renorm}
{\rm (Fonf)} Let $E$ be a separable Banach space.  Then $E$ is
isomorphically polyhedral if and only if $E'$ contains a countable
\ipn subset.
\end{thm}

This paper is devoted mainly to the problem of identifying the
isomorphically polyhedral Orlicz sequence spaces.  In \S 1, we prove
a characterization theorem for isomorphically polyhedral Banach spaces
having a shrinking basis.  This result is applied in \S 2 to obtain
examples of isomorphically polyhedral Orlicz spaces.  In \S 3, a
non-isomorphically polyhedral, $\cz$-saturated  Orlicz sequence space
is constructed.  Since every $\cz$-saturated Orlicz sequence space has
a separable dual, this shows that the converse of Theorem \ref{c0sat}
fails, answering a question posed by Rosenthal \cite{R1}.\\
\indent Standard Banach space terminology, as may be found in
\cite{LT}, is
employed. If $(e_n)$ is a basis of a Banach space $E$, and
$|||\cdot|||$ is a norm on $E$ equivalent to the given norm, we say
that $(e_n)$ is {\em monotone} with respect to $|||\cdot|||$ if
$|||\sum^{k}_{n=1}a_ne_n||| \leq |||\sum^{k+1}_{n=1}a_ne_n|||$ for every
real sequence $(a_n)$ and all $k \in \N$.  Terms and notation
regarding Orlicz spaces are discussed in \S 2.

\section{A characterization theorem}

This section is devoted to proving the following characterization
theorem.  Readers familiar with the proofs of Fonf's Theorems will
find the same ingredients used here.

\begin{thm}\label{char}
Let $(e_n)$ be a shrinking basis of a Banach space $(E, \|\cdot\|)$.
The following are equivalent.\\
(a) $E$ is isomorphically polyhedral;\\
(b) There exists an equivalent norm $|||\cdot|||$ on $E$ such that
$(e_n)$ is  a monotone basis with respect to $|||\cdot|||$, and for all\
$\sum a_ne_n \in E$, there exists $m \in {\rm \N}$ such that
\[ |||\sum^\infty_{n=1}a_ne_n||| = |||\sum^m_{n=1}a_ne_n|||. \]
\end{thm}

\begin{pf}
Let $(P_n)$ be the projections on $E$ associated with the basis
$(e_n)$.  The sequence $(P_n)$ is uniformly bounded with respect to
any equivalent norm on $E$.  Also, $(P_n)$ converges strongly to the
identity operator on $E$, which we denote by $1$.  Since $(e_n)$ is
shrinking, $(P'_n)$ converges to $1'$ strongly as well. \\
(a) $\Rightarrow$ (b).  By renorming, and using Theorem \ref{renorm},
we may assume that $E'$ contains a \pn sequence $(w_k)$.
Fix sequences $(\ep_k)$ and $(\delta_k)$ in $(0,1)$ which are both
convergent to $0$, and so that $(1 + \ep_k)(1 - 2\delta_k) > 1$ for
all $k$.  For each $k$, choose $n_k$ such that $\|(1 - P_n)'w_k\| \leq
\delta_k$ for all $n \geq n_k$.  Define a seminorm $|||\cdot|||$ on
$E$ by
\begin{equation}\label{def}
|||x||| = \sup_k(1+\ep_k)\max_{1\leq n\leq n_k}|\lan P_nx, w_k\ran|.
\end{equation}
Since $(w_k) \subseteq U_{E'}$, $|||x||| \leq 2\|x\|\sup\|P_n\|$.  On
the other hand, if $x \neq 0$, choose $k$ such that $\|x\| =
|w_k(x)|$.  Then
\begin{eqnarray*}
\|x\| & = & |w_k(x)| \hspace{1em} \leq \hspace{1em}
                          |\lan x, P'_{n_k}w_k\ran| +
                          |\lan x, (1 - P_{n_k})'w_k\ran| \\
& \leq & |\lan P_{n_k}x, w_k\ran| + \delta_k\|x\| .
\end{eqnarray*}
Thus
\begin{equation}\label{compare}
|||x||| \geq (1 + \ep_k)(1 - \delta_k)\|x\|  > \|x\| .
\end{equation}
Hence $|||\cdot|||$ is an equivalent norm on $E$.  It is clear that
$(e_n)$ is monotone with respect to $|||\cdot|||$.  We claim that this
norm satisfies the remaining condition in (b).  To this end, we first
show that the supremum in the definition (\ref{def}) is attained.
This is trivial if $x = 0$.  Fix
$0 \neq x \in E$.  Choose $k_1 \leq k_2 \leq \cdots$ and $(j_i)$ ,
$1 \leq j_i \leq n_{k_i}$ for all $i$, so that
\[ |||x||| = \lim_i(1 + \ep_{k_i})|\lan P_{j_i}x, w_{k_i}\ran|. \]
We divide the proof into cases. \\
\noindent\underline{Case 1}\hspace{1em} $\lim_ik_i = \lim_ij_i =
\infty$. \\
In this case, $P_{j_i}x \to x$ in norm.  Therefore
\[ \limsup_i|\lan P_{j_i}x, w_{k_i}\ran|  =
   \limsup_i|\lan x, w_{k_i}\ran| \leq \|x\| .\]
Also, $\ep_{k_i} \to 0$ as $i \to \infty$.  Thus, $|||x||| \leq
\|x\|$, contrary to (\ref{compare}).\\
\noindent\underline{Case 2}\hspace{1em} $\lim_ik_i = \infty$,
$\lim_ij_i \neq \infty$. \\
By using a subsequence, we may assume that $j_i = j$ for all $i$.  Then
\[ |||x||| = \lim_i(1 + \ep_{k_i})|\lan P_jx, w_{k_i}\ran|
      \leq \|P_jx\|. \]
Now choose $k$ such that $\|P_jx\| = |\lan P_jx, w_k\ran|$.  If $j
\leq n_k$,
\begin{eqnarray*}
|||x||| & \geq & (1 + \ep_k)|\lan P_jx, w_k\ran| \\
        & = & (1 + \ep_k)\|P_jx\| > \|P_jx\| ,
\end{eqnarray*}
a contradiction.  Now assume $j > n_k$, then
\begin{eqnarray*}
\|(P_j - P_{n_k})'w_k\| & \leq & \|(1 - P_j)'w_k\| + \|(1 -
P_{n_k})'w_k\| \\
& \leq & 2\delta_k .
\end{eqnarray*}
Hence
\begin{eqnarray*}
\|P_jx\| & = & |\lan P_jx, w_k\ran| \\
     & \leq & |\lan P_{n_k}x, w_k\ran| + 2\delta_k\|x\| \\
     & \leq & (1 + \ep_k)^{-1}|||x||| + 2\delta_k|||x|||.
\end{eqnarray*}
Therefore,
\[ |||x||| \leq \|P_jx\| \leq ((1 + \ep_k)^{-1} + 2\delta_k)|||x|||
    < |||x||| , \]
reaching yet another contradiction.  Consequently, we must have\\
\noindent\underline{Case 3}\hspace{1em} $\lim_ik_i \neq \infty$. \\
By using a subsequence, we may assume that the sequence $(k_i)$ is
constant.  Then it is clear that the supremum in (\ref{def})
is attained.   \\
\indent Now for any $x \in E$, choose $k$ so that the supremum in
(\ref{def})
is attained at $k$.  Then it is clear that $|||x||| =
|||P_{n_k}x|||$.\\

\noindent(b) $\Rightarrow$ (a).  Let $(\eta_n)$ and $(\ep_n)$ be
sequences convergent to $0$, with $1 > \eta_n > \ep_n > 0$ for all
$n$.  For each $n$, there is a finite $W_n \subseteq
U_{(E,|||\cdot|||)'}$ such that
\begin{equation}\label{Wn}
 (1 + \ep_n)^{-1}|||x||| \leq \max_{w\in W_n}|w(x)| \leq |||x|||
\end{equation}
for all $x \in \spn\{e_1, \ldots, e_n\}$.  Define a seminorm $\rho$ on
$E$ by
\begin{equation}\label{rho}
 \rho(x) = \sup_n(1 + \eta_n) \max_{1\leq j \leq n}\max_{w\in
                        W_j}|\lan P_jx, w\ran|.
\end{equation}
We will show that $\rho$ is an equivalent norm on $E$, and the set
\[W = \{(1 + \eta_n)P'_jw : n \in \N, 1 \leq j \leq n, w \in W_j\} \]
is a countable \pn subset of $(E, \rho)'$.  Then $E$ is isomorphically
polyhedral by Fonf's Theorem (Theorem \ref{renorm}).  Now let $x \in
E$. By (b), there exists $m$ such that $|||x||| = |||P_mx|||$.  Hence,
by (\ref{Wn}), and the fact that $(e_n)$ is monotone with respect to
$|||\cdot|||$,
\begin{eqnarray}\label{ineq}
|||x||| & = & |||P_mx||| \hspace{1em} \leq \hspace{1em}
              (1 + \ep_m)\max_{w\in W_m}|\lan P_mx, w\ran| \nonumber\\
        & \leq & (1 + \eta_m)\max_{w\in W_m}|\lan P_mx, w\ran| \\
        & \leq & \rho(x) \hspace{1em} \leq \hspace{1em} 2|||x|||.
                                   \nonumber
\end{eqnarray}
Thus $\rho$ is an equivalent norm on $E$. Next we show that the
supremum in (\ref{rho}) is attained. Fix $x \in E$.  Choose sequences
$n_1 \leq n_2 \leq \cdots$, $(j_k)$, and $(w_k)$ such that $1 \leq j_k
\leq n_k$, $w_k \in W_{n_k}$ for all $k$, and
$\rho(x) = \lim_k(1 + \eta_{n_k})|\lan P_{j_k}x,w_k\ran|$.
First assume that $\lim_kn_k = \infty$.  Then $\eta_{n_k}
\to 0$.  Since $(e_n)$ is monotone with respect to $|||\cdot|||$, we
have $\rho(x) \leq |||x|||$.  But there exists $k$ such that $|||x|||
= |||P_kx|||$, and there is a $w \in W_k$ such that $|||P_kx||| \leq
(1 + \ep_k)|w(P_kx)|.$  Thus
\[ \rho(x) \geq (1 + \eta_k)|w(P_kx)| \geq \frac{1 + \eta_k}{1 +
\ep_k}|||x||| > |||x|||, \]
a contradiction.  Therefore, $\lim_kn_k \neq \infty$.  By going to a
subsequence, we may assume that $(n_k)$ is
bounded. Using a further subequence if necessary, we may even assume
it is constant.  Thus the supremum in (\ref{rho}) is attained.  From
this it readily follows that the set $W$ is a \pn subset of $(E,
\rho)'$.  The countability of $W$ is evident.
\end{pf}

\noindent {\bf Remark} The assumption that the basis $(e_n)$ is
shrinking is used only in the proof of (a) $\Rightarrow$ (b).  If
$(e_n)$ is assumed to be unconditional and (a) holds, then $(e_n)$
must be shrinking.  For otherwise $E$ contains a copy of $\ell^1$,
which contradicts (a) by Fonf's Theorem (Theorem \ref{c0sat}).
  Thus the assumption of
shrinking is not needed if $(e_n)$ is unconditional.\\

\section{Orlicz sequence spaces}

In this section, we apply Theorem \ref{char} to identify a class of \iso
\poly Orlicz sequence spaces.  Terms and notation about Orlicz
sequence spaces follow that of \cite{LT}.  An {\em Orlicz function}
$M$ is
a continuous non-decreasing convex function defined for $t \geq 0$
such that $M(0) = 0$ and $\lim_{t\to \infty}M(t) = \infty$.  If $M(t)
> 0$ for all $t > 0$, then it is {\em non-degenerate}.  Clearly a
non-degenerate Orlicz function must be strictly increasing. The {\em
Orlicz sequence space} $\ell_M$ associated with an Orlicz function $M$
is the space of all sequences $(a_n)$ such that $\sum M(|a_n|/\rho)
< \infty$ for some $\rho > 0$, equipped with the norm
\[ \|x\| = \inf\{\rho > 0 : \sum M(|a_n|/\rho) < \infty\} .\]
Let $e_n$ denote the vector whose sole nonzero coordinate is a $1$ at
the $n$-th position.  Then clearly $(e_n)$ is a basic sequence in
$\ell_M$. The closed linear span of $\{e_n\}$ in $\ell_M$ is denoted
by $h_M$. Alternatively, $h_M$ may be described as the set of all
sequences $(a_n)$ such that $\sum M(|a_n|/\rho) < \infty$ for every
$\rho > 0$.  Additional results and references on Orlicz spaces may be
found in \cite{LT}. For a real null sequence $(a_n)$, let $(a^*_n)$
denote the decreasing rearrangement of the sequence $(|a_n|)$.

\begin{thm}\label{ip}
Let $M$ be a non-degenerate Orlicz function such that there exists a
finite number K satisfying $\lim_{t\to 0}M(Kt)/M(t) = \infty$.  Then
$h_M$ is isomorphically polyhedral.
\end{thm}

\begin{pf}
For all $k \in \N$, let
\[ b_k = \inf\llb\frac{M(Kt)}{M(t)} : 0 < t \leq
M^{-1}(\frac{1}{k})\lrb.
\]
Then $\lim_{k\to\infty} b_k = \infty$.  Thus there is a sequence
$(\eta_k)$ decreasing to $1$ such that $\eta_k > (1 -
b^{-1}_{k+1})^{-1}$ for all $k$.  Define a seminorm on $h_M$ by
\begin{equation}\label{defO}
|||(a_n)||| = \sup_k\eta_k\|(a^*_1, \ldots, a^*_k, 0, \ldots)\| ,
\end{equation}
where $\|\cdot\|$ is the given norm on $h_M$.  It is clear that
$|||\cdot|||$ is an equivalent norm on $h_M$, and that $(e_n)$ is a
monotone basis with respect to $|||\cdot|||$. It suffices to show that
$|||\cdot|||$ satisfy the remaining condition in part (b) of Theorem
\ref{char}.  We first show that if $(a_n)$ is a positive decreasing
sequence in $h_M$, then there is a $k$ such that
\begin{equation}\label{growth}
\|(a_n)\| \leq \eta_k\|(a_1, \ldots, a_k, 0, \ldots)\| .
\end{equation}
Assume otherwise.  There is no loss of generality in assuming that
$\|(a_n)\| = 1$.  Then $\sum M(a_n) = 1$ and $\sum^k_{n=1}M(\eta_ka_n)
\leq 1$ for all $k$.  In particular, note that the second condition
implies $a_k \leq M^{-1}(1/k)$ for all $k$, since $\eta_k \geq 1$ and
$(a_n)$ is decreasing.
Now choose $m$ such that $\|(0, \ldots, 0, a_m,
a_{m+1}, \ldots)\| \leq K^{-1}$.  Then $\summ M(Ka_n) \leq 1$.  Also
$M(Ka_n) \geq b_mM(a_n)$ for all $n \geq m$.  Therefore,
\begin{eqnarray*}
1 & = & \sum M(a_n) \hspace{1em} = \hspace{1em} \sumum M(a_n) + \summ
M(a_n) \\
& \leq & \eta^{-1}_{m-1}\sumum M(\eta_{m-1}a_n) + b^{-1}_m\summ
M(Ka_n) \\
& \leq & \eta^{-1}_{m-1} + b^{-1}_m \hspace{1em} < \hspace{1em} 1 ,
\end{eqnarray*}
a contradiction. Hence (\ref{growth}) holds for some $k$.
Now for a general element $(a_n) \in h_M$, choose $m$ such that
$\|(a_n)\| = \|(a^*_n)\| \leq \eta_m\|(a^*_1, \ldots, a^*_m, 0,
\ldots)\|$.    Note that since $\lim_k\eta_k\|(a^*_1, \ldots, a^*_k, 0,
\ldots)\| = \|(a_n)\|$, the supremum in equation (\ref{defO}) is
attained, say, at $j$.  Then choose $i$ large enough that $a^*_1,
\ldots, a^*_j$ are found in $\{|a_1|, \ldots, |a_i|\}$.  With this
choice of $i$,
\[ |||(a_1, \ldots, a_i, 0, \ldots)||| \geq
\eta_j\|(a^*_1, \ldots, a^*_j, 0,
\ldots)\| = |||(a_n)|||
\]
by choice of $j$.  Since the reverse inequality is obvious,
\[ |||(a_n)||| = |||(a_1, \ldots, a_i, 0, \ldots)|||, \]
as required.
\end{pf}

\section{A counterexample}

\begin{thm}\label{ex}
Let $M$ be a non-degenerate Orlicz function.  Suppose there exists a
sequence $(t_n)$ decreasing to $0$ such that
\[ \sup_n \frac{M(Kt_n)}{M(t_n)} < \infty \]
for all $K < \infty$.  Then $h_M$ is not isomorphically polyhedral.
\end{thm}

\begin{pf}
Suppose that $h_M$ is \iso \poly.  By Theorem \ref{char} and the
remark following it, one obtains a norm $|||\cdot|||$ on $h_M$ as
prescribed by part (b) of the theorem.  Fix $\alpha > 0$ so that
$|||x||| \leq \alpha \Rightarrow \|x\| \leq 1$.  Choose a sequence
$(\eta_k)$ strictly decreasing to $1$.  Let $n_1 = \min\{n\in \N :
\eta_1|||t_ne_1||| \leq \alpha\}$.  If $n_1 \leq n_2 \leq \cdots \leq
n_k$  are
chosen so that $\eta_k|||\sumnj||| \leq \alpha$, then
$\eta_{k+1}|||\sumnj||| < \alpha$.  Hence
\[ \{n \geq n_k : \eta_{k+1}|||\sumnj + t_ne_{k+1}||| \leq \alpha\}
\neq \emptyset . \]
Now define
\begin{equation}\label{min}
n_{k+1} = \min\{n \geq n_k : \eta_{k+1}|||\sumnj + t_ne_{k+1}||| \leq
\alpha\} .
\end{equation}
This inductively defines a (not necessarily strictly) increasing
sequence $(n_k)$ satisfying
\begin{equation}\label{eta}
\eta_k|||\sumnj||| \leq \alpha
\end{equation}
for all $k$ and  the minimality condition (\ref{min}).
In particular, $|||\sumnj||| \leq \alpha$ for all $k$, so $\|\sumnj\|
\leq 1$ by
the choice of $\alpha$. Therefore $\sumMnj \leq 1$ for all $k$.  For
all $K < \infty$ and all $k \in \N$,
\[ \sumj M(Kt_{n_j}) \leq \sup_m\frac{M(Kt_m)}{M(t_m)}\sumj M(t_{n_j})
\leq \sup_m\frac{M(Kt_m)}{M(t_m)} .\]
Consequently, $\sum^\infty_{j=1}M(Kt_{n_j}) < \infty$ for all $K <
\infty$.   Hence $x = \sum^\infty_{j=1}t_{n_j}e_j$ converges in $h_M$.
Clearly $|||x||| = \lim_k|||\sumnj||| \leq \alpha$.  We claim that in
fact $|||x||| = \alpha$.  Otherwise, suppose $|||x||| = \beta <
\alpha$. Since $(e_n)$ is monotone with respect to $|||\cdot|||$,
$|||\sumnj||| \leq \beta < \alpha$ for all $k$.  By the convergence of
$x$, $\lim_jt_{n_j} = 0$.  So one can find $i$ such that
$|||t_{n_i}e_j||| \leq \alpha - \beta$ for all $j$.  Then
\[ |||\sumi + t_{n_i}e_{i+1}||| \leq |||\sumi||| +
|||t_{n_i}e_{i+1}||| \leq \beta + \alpha - \beta = \alpha .
\]
By the minimality condition (\ref{min}), $n_{i+1} = n_i$.  Similarly,
we see that $n_j = n_i$ for all $j \geq i$.  This contradicts the
convergence of $x$ and proves the claim.  But now, by
(\ref{eta}), $|||\sumnj||| < \alpha = |||x|||$ for all $k$,
contradicting the choice of the norm $|||\cdot|||$.
\end{pf}

We now construct an Orlicz function $M$ satisfying Theorem \ref{ex}
while $h_M$ is $\cz$-saturated.  We begin with some simple results
which help to identify the $\cz$-saturated Orlicz sequence spaces.

\begin{pr}
Let $M$ be a non-degenerate Orlicz function.  Then the following are
equivalent.\\
(a) $h_M$ is $\cz$-saturated;\\
(b) $h_M$ does not contain an isomorph of $\ell^p$ for any $1 \leq p <
\infty$; \\
(c) for all $q < \infty$,
\[\sup_{0<\lambda,t\leq 1}\frac{M(\lambda t)}{M(\lambda)t^q} < \infty.
\]
\end{pr}

\begin{pf}
Clearly (a) implies (b).  If (a) fails, let $Y$ be an infinite
dimensional closed subspace of $h_M$ which contains no isomorph of
$\cz$. By \cite[Proposition 4.a.7]{LT}, $Y$ has a subspace $Z$
isomorphic to some Orlicz sequence space $h_N$.  Then $h_N$ contains
no isomorph of $\cz$.  By \cite[Theorem 4.a.9]{LT}, $h_N$ contains an
isomorph of some $\ell^p, 1 \leq p < \infty$.  Hence $Y$ contains a
copy of $\ell^p$, and (b) fails.  The equivalence of (b) and (c)
also follows from \cite[Theorem 4.a.9]{LT}.
\end{pf}

\begin{pr}\label{bees}
Let $(b_n)^\infty_{n=0}$ be a decreasing sequence of strictly positive
numbers such that
\[\sup_{m,n}\frac{b_{m+n}}{b_n}K^m < \infty \hspace{1em} \mbox{for
all} \hspace{1em} K < \infty. \]
Define $M$ to be the continuous, piecewise linear function such that
$M(0) = 0$,
\[ M'(t) = \left\{ \begin{array}{ll}
                    b_n & \mbox{if\hspace{.8em}$2^{-n-1} < t <
                                      2^{-n}$,  $n > 0$}\\
                    b_0 & \mbox{if\hspace{.8em}$ 2^{-1} < t$}
                   \end{array}
           \right. \]
Then $M$ is a non-degenerate Orlicz function so that $h_M$ is
$\cz$-saturated.
\end{pr}

\begin{pf}
It is clear that $M$ is a non-degenerate Orlicz function.  For all $n
\geq 0$, $2^{-n-1}b_n \leq M(2^{-n}) \leq 2^{-n}b_n$.  Hence
\[ C_q \equiv \sup_{m,n}\frac{M(2^{-m-n})}{M(2^{-n})}2^{mq}
\leq 2\sup_{m,n}\frac{b_{m+n}}{b_n}(2^{q-1})^m < \infty \]
for any $q < \infty$.  Now if $\lambda, t \in (0,1]$, choose $m, n
\geq 1$ such that $t \in (2^{-m}, 2^{-m+1}]$, $\lambda \in (2^{-n},
2^{-n+1}]$.  Then $\lambda t \in (2^{-m-n}, 2^{-m-n+2}]$.  If $m \geq
2$, then
\[ \frac{M(\lambda t)}{M(\lambda)t^q} \leq
2^{2q}\frac{M(2^{-(m-2)-n})}{M(2^{-n})}2^{(m-2)q} \leq 4^qC_q.\]
If $m = 1$, then $t > 2^{-1}$.  Therefore
\[ \frac{M(\lambda t)}{M(\lambda)t^q} \leq t^{-q} \leq 2^q. \]
Thus
\[\sup_{0<\lambda,t\leq 1}\frac{M(\lambda t)}{M(\lambda)t^q} < \infty,
\]
and $h_M$ is $\cz$-saturated by the previous proposition.
\end{pf}

\begin{thm}
There exists an Orlicz function $M$ such that $h_M$ is $c_0$-saturated
but not isomorphically polyhedral.  In particular, a $c_0$-saturated
space with a separable dual is not necessarily isomorphically
polyhedral.
\end{thm}

\begin{pf}
It is well known that every $\cz$-saturated space $h_M$ has a
separable dual.  Thus the second statement follows from the first.
Let $\alpha_0 = \alpha_1 = \alpha_2 = 1$, and let $\alpha_j = (e/j)^j$
for $j \geq 3$.  Then $(\alpha_j)$ is a decreasing sequence.  Choose a
decreasing sequence $(c_j)^\infty_{j=0}$ of strictly positive numbers
such that $c_{j+1} \leq \alpha_j\alpha_{2j^2}c_j$ for all $j \geq 0$.
For convenience, set $s_n = \sum^n_{j=1}j$ for all $n \geq 1$.  Now
define $b_0 = c_0$, $b_1 = c_1$, and $b_{s_n+k} =
c_{n+1}/\alpha_{n+1-k}$ whenever $n \geq 1$ and $1 \leq k \leq n+1$.
We first show that the sequence $(b_j)$ satisfies the conditions in
Proposition \ref{bees}.  \\
\underline{Claim 1}\hspace{1em}$(b_j)$ is a decreasing sequence.  \\
One verifies directly that $b_0 \geq b_1 \geq b_2$.  If $n \geq 1$ and
$1 \leq k \leq j \leq n+1$,
\[b_{s_n+k} = \frac{c_{n+1}}{\alpha_{n+1-k}} \geq
\frac{c_{n+1}}{\alpha_{n+1-j}} = b_{s_n+j} \]
since $(\alpha_m)$ is decreasing. Finally,
\[ b_{s_{n+1}+1} = \frac{c_{n+2}}{\alpha_{n+1}} \leq
\alpha_{2(n+1)^2}c_{n+1} \leq c_{n+1} = b_{s_n+n+1} \]
for all $n \geq 1$.  This proves Claim 1.\\
\underline{Claim 2}\hspace{1em} $b_{m+n} \leq \alpha_mb_n$ for all $m
\geq 0$, $n \geq 2$.\\
Express $n = s_i + k$, $m + n = s_j + l$, where $1 \leq i \leq j$, $1
\leq k \leq i+1$, and $1 \leq l \leq j + 1$.  If $i = j$, then $l -
k = m$. Moreover, $i + 1
- k \geq \max\{l - k, i + 1 - l\}$, from which it follows
that $\alpha_{i+1-k} \leq \alpha_{l-k}\alpha_{i+1-l}$. Therefore,
\[ b_{m+n} = \frac{c_{i+1}}{\alpha_{i+1-l}} \leq
             \alpha_m\frac{c_{i+1}}{\alpha_{i+1-k}} = \alpha_mb_n .\]
Now consider the possibility that $j > i$. Note first that
\[ m = (m + n ) - n \leq s_j + j + 1 - (s_i + 1) \leq s_j + j \leq
2j^2. \]
Hence $\alpha_m \geq \alpha_{2j^2}$.  Using Claim 1 and the properties
of the sequence $(c_j)$, we obtain
\begin{eqnarray*}
b_{m+n} & \leq & b_{s_j+1} \hspace{1em} = \hspace{1em}
                 \frac{c_{j+1}}{\alpha_j} \\
  & \leq & \alpha_{2j^2}c_j \hspace{1em} \leq \hspace{1em}
                 \alpha_mc_{i+1} \\
  & = & \alpha_mb_{s_i+i+1} \fatleq \alpha_mb_n.
\end{eqnarray*}
\noindent\underline{Claim 3}\hspace{1em}
\[\sup_{m,n}\frac{b_{m+n}}{b_n}K^m < \infty \hspace{1em} \mbox{for
all} \hspace{1em} K < \infty. \]
First observe that for $i \geq 1$, $1 \leq k \leq i + 1$, and $K <
\infty$,
\begin{eqnarray*}
 b_{s_i+k}K^{s_i+k} & = & \frac{c_{i+1}}{\alpha_{i+1-k}}K^{s_i+k}
              \hspace{1em}
        \leq \hspace{1em} \alpha_{2i^2}c_iK^{s_i+i+1} \\
         & \leq & c_0\alpha_{2i^2}K^{s_i+i+1} \fatarr 0
\end{eqnarray*}
as $i \to \infty$.  Hence $(b_mK^m)_m$ is bounded.  Therefore
$\sup_{n=1,2}\sup_{m}b_{m+n}K^m/b_n < \infty$.  On the other hand,
using Claim 2,
\[ \sup_{n\geq 2}\sup_{m}\frac{b_{m+n}}{b_n}K^m \leq
       \sup_{m}\alpha_mK^m < \infty \]
by direct verification.\\
\indent Define the function $M$ using the sequence $(b_j)$ as in
Proposition \ref{bees}.  Using Claims 1 and 3, and the proposition, we
see that $h_M$ is $\cz$-saturated.  To complete the proof, it suffices
to find a
sequence $(t_n)$ as in Theorem \ref{ex}.  We claim that the sequence
$(t_n) = (2^{-s_n})$ will do.  Clearly $(t_n)$ decrease to $0$.  Fix
$m \in \N$.  For all $n > m$,
\[ b_{s_n-m} = b_{s_{n-1}+(n-m)} = \frac{c_n}{\alpha_m}. \]
Hence
\begin{eqnarray*}
M(2^mt_n) & = & M(2^{-s_n+m}) \fatleq \frac{b_{s_n-m}}{2^{s_n-m}} \\
 & = & \frac{c_n}{\alpha_m2^{s_n-m}} \fateq
     \frac{2^{m+1}}{\alpha_m}\frac{c_n}{2^{s_n+1}} \\
 & = & \frac{2^{m+1}}{\alpha_m}\frac{b_{s_n}}{2^{s_n+1}} \fatleq
      \frac{2^{m+1}}{\alpha_m}M(t_n)
\end{eqnarray*}
whenever $n > m$.  Therefore,
\[ \sup_n \frac{M(2^mt_n)}{M(t_n)} < \infty \]
for all $m \in \N$.
\end{pf}

The obvious question to be raised is how to characterize \iso \poly
$h_M$ in terms of the Orlicz function $M$.  We suspect that the
condition given in Theorem \ref{ip} is the correct one.  It can be
shown that if $\liminf_{t\to 0} M(Kt)/M(t) < \infty$ for all $K <
\infty$, then for any sequence $(\eta_k)$ decreasing to $1$, the norm
given by equation (\ref{defO}) does not satisfy part (b) of Theorem
\ref{char}.

\flushleft
\vspace{.5in}
Department of Mathematics\\National University of Singapore\\
Singapore 0511\\ e-mail(bitnet) : matlhh@nusvm

\end{document}